# On the frequentist coverage of Bayesian credible intervals for lower bounded means

### Éric Marchand

*Département de mathématiques, Université de Sherbrooke, Sherbrooke, Qc, CANADA, J1K 2R1*
*e-mail:* eric.marchand@usherbrooke.ca

### William E. Strawderman

*Department of Statistics, Rutgers University, 501 Hill Center, Busch Campus, Piscataway, N.J., USA 08854-8019*
*e-mail:* straw@stat.rutgers.edu

### Keven Bosa

*Statistics Canada - Statistique Canada, 100, promenade Tunney's Pasture, Ottawa, On, CANADA, K1A 0T6*
*e-mail:* keven.bosa@usherbrooke.ca

### Aziz Lmoudden

*Département de mathématiques, Université de Sherbrooke, Sherbrooke, Qc, CANADA, J1K 2R1*
*e-mail:* aziz.lmoudden@usherbrooke.ca

**Abstract:** For estimating a lower bounded location or mean parameter for a symmetric and logconcave density, we investigate the frequentist performance of the $100(1-\alpha)\%$ Bayesian HPD credible set associated with priors which are truncations of flat priors onto the restricted parameter space. Various new properties are obtained. Namely, we identify precisely where the minimum coverage is obtained and we show that this minimum coverage is bounded between $1 - \frac{3\alpha}{2}$ and $1 - \frac{3\alpha}{2} + \frac{\alpha^2}{1+\alpha}$; with the lower bound $1 - \frac{3\alpha}{2}$ improving (for $\alpha \leq 1/3$) on the previously established ([9]; [8]) lower bound $\frac{1-\alpha}{1+\alpha}$. Several illustrative examples are given.

## Contents

É. Marchand et al./Frequentist coverage of Bayesian credible intervals 10293   Examples and final comments . . . . . . . . . . . . . . . . . . . . . . 1038
4   Appendix . . . . . . . . . . . . . . . . . . . . . . . . . . . . . . . . . 1041
Acknowledgements . . . . . . . . . . . . . . . . . . . . . . . . . . . . . . 1042
References . . . . . . . . . . . . . . . . . . . . . . . . . . . . . . . . . . 1042
## 1. Introduction

The findings in this paper are concerned with interval estimation of a location parameter $\theta$ based on an observable $X$ having cdf $G(x - \theta)$ with $\theta \geq 0$ and symmetric densities $g(x - \theta)$, in other words for cases where there exists a lower bound constraint for $\theta$ with the parameter space's minimal value taken to be 0 without loss of generality. More specifically, we establish new analytical properties of the frequentist coverage

$$C(\theta) = P_\theta(I_{\pi^*}(X) \ni \theta), \ (\theta \geq 0), \tag{1}$$

of the $100 \times (1-\alpha)\%$ Bayesian HPD confidence interval $I_{\pi^*}(X)$ associated with the flat prior $\pi^*(\theta) = 1_{[0,\infty)}(\theta)$.

With a known but conservative $\frac{1-\alpha}{1+\alpha}$ lower bound for the minimal coverage $\inf_{\theta \geq 0} C(\theta)$ (see [9] for the normal case; [8] for the general unimodal case), we succeed here in establishing for logconcave pdf's $g$ the better (for $\alpha < 1/3$) lower bound $1 - \frac{3\alpha}{2}$ for $\inf_{\theta \geq 0} C(\theta)$. We further establish various new properties of the frequentist coverage $C(\theta)$ for logconcave pdf's $g$, such as **(i)** $argmin_{\theta \geq 0} C(\theta) = 2d_0$; **(ii)** $1 - \frac{3\alpha}{2} \leq \inf_{\theta \geq 0} C(\theta) \leq 1 - \frac{3\alpha}{2} + \frac{\alpha^2}{1+\alpha}$; **(iii)** $\sup_{\theta \geq 0} C(\theta) \leq 1 - \frac{\alpha}{2}$; and **(iv)** $C(\theta)$ decreases on $(d_1, 2d_0)$, increases on $(2d_0, \infty)$, where $d_0 = G^{-1}(\frac{1}{1+\alpha})$ and $d_1 = G^{-1}(1 - \frac{\alpha}{2})$. Observe that result **(ii)** confirms the first order optimality of the lower bound $1 - \frac{3\alpha}{2}$, as opposed for instance to the lower bound $\frac{1-\alpha}{1+\alpha} = 1 - 2\alpha + o(\alpha)$ for small $\alpha$. The main application of our results is undoubtedly for estimating a lower bounded normal mean with known variance, but the results are nevertheless more general and unified for symmetric densities with logconcave pdf's. For instance, an interesting illustration will be given for a Laplace model.

Here is a glimpse of previous work pertaining to the credible interval $I_{\pi^*}(X)$ and its frequentist coverage $C(\theta)$. For estimating a lower bounded normal mean $\theta$ (say $\theta \geq 0$), results due to **(a)** Roe and Woodroofe [9] (known variance), and **(b)** Zhang and Woodroofe [10] (unknown variance $\sigma^2$) establish the lower bound $\frac{1-\alpha}{1+\alpha}$ for the frequentist coverage of the $100 \times (1-\alpha)\%$ HPD Bayesian confidence interval with respect to the priors: $\pi(\theta) = 1_{[0,\infty)}(\theta)$ in **(a)**; and $\pi_0(\theta, \sigma) = \frac{1}{\sigma} 1_{[0,\infty)}(\theta) 1_{(0,\infty)}(\sigma)$ in **(b)**. More recently, Marchand and Strawderman [8] established, for a more general setting with underlying symmetry, the validity of the lower bound $\frac{1-\alpha}{1+\alpha}$ for the frequentist coverage of the $100 \times (1-\alpha)\%$ Bayesian HPD credible interval derived from the truncation onto the restricted parameter space of the Haar right invariant prior. We refer to their paper for details and similar developments for non-symmetric settings (also see [11]; [12]). The important starting point to keep in mind from the results of Marchand and Strawderman [8] is the applicability of the $\frac{1-\alpha}{1+\alpha}$ lower bound for symmetric and unimodal location models.

As mentioned above, the analysis which we have carried out and which is presented herein was motivated by the conservative nature of the previously established lower bound $\frac{1-\alpha}{1+\alpha}$ and other unknown and potentially useful aspects of the coverage $C(\theta)$. Finally, the results are cast here with the backdrop of a flurry of recent activity and debate, which is focused on the choice of methods for setting confidence bounds for restricted parameters as witnessed by the works referred to above as well at that of Mandelkern [7], Feldman and Cousins [5], Efron [4], among others.

## 2. Main results

Our main results apply to location models $X \sim g(x - \theta)$ with densities $g$ which are unimodal, symmetric about $0$, and logconcave (i.e., $\log g$ is a concave function on its support). For such location families, the key assumption of logconcavity can equivalently be described as corresponding to those families with increasing monotone likelihood ratio densities, that is $\frac{g(x-\theta_1)}{g(x-\theta_0)}$ increases in $x$ for all $\theta_1, \theta_0$ such that $\theta_1 > \theta_0$. Alternatively, the assumption of logconcavity is connected to an increasing hazard (or failure) rate (see Lemma 1). Before moving along, here is a useful checklist of notations and definitions used throughout.

*Checklist*

- $1 - \alpha$: credibility or posterior coverage or nominal frequentist coverage ($\alpha \in (0, 1)$)
- $g$: probability density function (pdf) of $X - \theta$
- $G$: cumulative distribution function (cdf) of $X - \theta$
- $\bar{G}(\equiv 1 - G)$: survivor function
- $G^{-1}$: inverse cdf
- $\lambda(\cdot)$: hazard rate function given by $\lambda(z) = \frac{g(z)}{\bar{G}(z)}$
- $\pi^*$: the flat prior density truncated onto the parameter space $[0, \infty)$ given by $\pi^*(\theta) = 1_{[0,\infty)}(\theta)$
- $I_{\pi^*}(X) = [l(X), u(X)]$: the HPD Bayesian credible set of credibility $1 - \alpha$ associated with $\pi^*$
- $a = \lim_{x \to -\infty} u(x)$
- $C(\theta)$: the frequentist coverage at $\theta$ of the confidence interval $I_{\pi^*}(X)$ given by $C(\theta) = P_\theta(I_{\pi^*}(X) \ni \theta)$
- $d_0 = G^{-1}(\frac{1}{1+\alpha})$
- $d_1 = G^{-1}(1 - \frac{\alpha}{2})$
- $u^{-1}(\theta) = \inf\{x|\theta \leq u(x)\}$
- $l^{-1}(\theta) = \sup\{x|\theta \geq l(x)\}$

(these last two inverses being defined as such in order to facilitate expressions below for frequentist coverage)

We next collect some useful properties of logconcave densities (see [1], [2] for surveys). Lemma 2 is a critical result which may well be new and of independent interest.

**Lemma 1.** *Let $g$ be a unimodal, symmetric about $0$ and logconcave density on $\mathbb{R}$. Then,*

(a) $\frac{g(y)}{g(y+\theta)}$ *is nondecreasing in $\theta$ on $[0,\infty)$, for $y \geq 0$;*
(b) $\frac{g(y)}{g(y+\theta)}$ *is nondecreasing in $y$ on $\mathbb{R}$, for $\theta \geq 0$;*
(c) $\frac{g(y(\theta))}{g(y(\theta)+\theta)}$ *is nondecreasing in $\theta$ on $[0,\infty)$, for $y(\cdot)$ a nonnegative and nondecreasing function taking values on $[0,\infty)$;*
(d) *both the cdf $G$ and survivor function $\overline{G}$ are log-concave on $[0,\infty)$;*
(e) *the hazard rate $\lambda(\cdot)$ is nondecreasing on $\mathbb{R}$.*

*Proof.* Part (a) is obvious and is simply a consequence of unimodality. Parts (b) and (e) are log-concavity, and part (c) follows from (a) and (b). A proof of (d) can be found in [1] or [2]. □

**Lemma 2.** *Let $g$ be a unimodal, symmetric about $0$ and logconcave density on $\mathbb{R}$. Then, for all $z \geq 0$,*

(a)
$$\frac{g(z)}{g(2z)} \geq \frac{1}{2\overline{G}(z)} - 1\,; \qquad (2)$$

(b)
$$\frac{\overline{G}(z)}{\overline{G}(2z)} \geq \frac{1}{2\overline{G}(z)} - 1\,. \qquad (3)$$

*Proof.* Part (a) implies (b) given the increasing hazard rate property of part (e) of Lemma 1, as

$$\lambda(2z) \geq \lambda(z) \implies \frac{\overline{G}(z)}{\overline{G}(2z)} \geq \frac{g(z)}{g(2z)} \text{ for all } z \geq 0\,.$$

For part (a), begin by observing that for all $z > 0$:

$$\frac{1}{2\overline{G}(z)} - 1 \leq \frac{\frac{1}{2} - \overline{G}(z)}{\overline{G}(z) - \overline{G}(2z)} = \frac{\int_0^z g(y)\,dy}{\int_z^{2z} g(y)\,dy},$$

so that (2) holds as soon as, for all $z > 0$,

$$\int_0^z \frac{g(y)}{g(z)}\,dy \leq \int_z^{2z} \frac{g(y)}{g(2z)}\,dy,$$

or equivalently,

$$\int_0^z \left(\frac{g(y)}{g(z)} - \frac{g(y+z)}{g(2z)}\right) dy \leq 0. \qquad (4)$$

Finally, apply part (b) of Lemma 1 inside the above integral to infer that $\frac{g(y)}{g(z)} - \frac{g(y+z)}{g(2z)} \leq 0$ for all $y, z$ such that $y \in (0, z)$, which yields (4) and the proof of (2). □

We now recall previously established properties of $I_{\pi^*}(X)$ and of its frequentist coverage before pursuing with the main analysis.

**Lemma 3** ([8], theorem 1). *For $X \sim g(x - \theta)$, $\theta \geq 0$, $g$ unimodal, symmetric about 0, and the HPD credible set $I_{\pi^*}(X)$, we have*

(a) $I_{\pi^*}(X) = [l(X), u(X)]$, with $l(x) = \{x - G^{-1}(\frac{1}{2} + \frac{1-\alpha}{2}G(x))\}1_{(d_0, \infty)}(x)$ and $u(x) = \{x - G^{-1}(\alpha G(x))\}1_{(-\infty, d_0]}(x) + \{x + G^{-1}(\frac{1}{2} + \frac{1-\alpha}{2}G(x))\}1_{(d_0, \infty)}(x)$;
(b) $C(\theta) > \frac{1-\alpha}{1+\alpha}$, for all $\theta \geq 0$;
(c) $C(0) = \frac{1}{1+\alpha}$;
(d) $\lim_{\theta \to \infty} C(\theta) = 1 - \alpha$.

As shown in the following lemma, the lower bound $l(\cdot)$ is nondecreasing on $(-\infty, \infty)$, while the upper bound $u(\cdot)$ is nondecreasing on $(d_0, \infty)$. Furthermore, we show how the logconcavity of $g$ forces $u(\cdot)$ to be nondecreasing on $(-\infty, d_0)$ as well. (It is easily verified that $l(\cdot)$ and $u(\cdot)$ are continuous functions on $\mathbb{R}$.)

**Lemma 4.** *Consider $I_{\pi^*}(X)$ as given in Lemma 3. Then,*

(a) $l(\cdot)$ is nondecreasing on $(-\infty, \infty)$;
(b) $u(\cdot)$ is nondecreasing on $(d_0, \infty)$;
(c) $u(\cdot)$ is nondecreasing on $(-\infty, d_0)$ whenever $g$ is logconcave.

*Proof.* **(a)** Since $l(x) = 0$ for $x \leq d_0$, we only need to look at the behaviour of $l(x)$ for $x > d_0 (> 0)$. We have, for $x > d_0$,

$$\frac{d}{dx} l(x) = 1 - \frac{1-\alpha}{2} \frac{g(x)}{g(G^{-1}(\frac{1}{2} + \frac{1-\alpha}{2}G(x)))}.$$

Notice, since $l(x) \geq 0$, that we must have $x \geq G^{-1}(\frac{1}{2} + \frac{1-\alpha}{2}G(x)) \geq 0$, from which we infer that $g(x) \leq g(G^{-1}(\frac{1}{2} + \frac{1-\alpha}{2}G(x)))$ given that $g$ is unimodal with a mode at 0.

**(b)** Follows directly with $G$ and $G^{-1}$ being nondecreasing.[1]

**(c)** Under the prior $\pi^*$, the posterior density of $\theta|x$ is given by $\{\frac{g(\theta-x)}{G(x)}1_{[0,\infty)}(\theta)\}$, and $u(x)$ is, for $x \leq d_0$, the corresponding quantile of order $1 - \alpha$ ([8], proofs of Theorem 1 or Lemma 5). Now the logconcavity of $g$ implies that the family of posterior densities of $\theta|x$, with parameter $x$, has increasing monotone likelihood ratio in $\theta$. Finally, the result follows since the quantiles $u(x); x \leq d_0$; of such families are nondecreasing in $x$.

The increasing properties of $l(\cdot)$ and $u(\cdot)$, for logconcave densities $g$, permit us to translate the events of coverage, i.e., $\{I_{\pi^*}(X) \ni \theta\}$, as $\{u^{-1}(\theta) \leq X \leq l^{-1}(\theta)\}$, and leads to the following expressions and further properties of the frequentist coverage $C(\theta)$. □

---
[1] actually, we have the stronger result that $u(x) - x$ increases in $x$ here for $x \geq d_0$

**Lemma 5.** *Under the conditions of Lemma 3 with g logconcave, we have:*

**(a)**

$$C(\theta) = \begin{cases} G(x_1(\theta)) & \text{if } \theta \in [0, a] \\ G(x_1(\theta)) - G(x_0(\theta)) & \text{if } \theta \in (a, 2d_0] \\ G(x_1(\theta)) - G(x_2(\theta)) & \text{if } \theta \in (2d_0, \infty) \end{cases} \quad (5)$$

*where $a = \lim_{x \to -\infty} u(x)$ (see checklist), and with $x_0(\cdot)$, $x_1(\cdot)$, and $x_2(\cdot)$ being functions defined by the equations:*

$$\begin{align} (i) \ & G(x_0(\theta)) = \alpha\, G(x_0(\theta) + \theta) & (\theta \in (a, 2d_0]); \\ (ii) \ & 2G(x_1(\theta)) - 1 = (1 - \alpha)G(x_1(\theta) + \theta) & (\theta \geq 0); \\ (iii) \ & 1 - 2G(x_2(\theta)) = (1 - \alpha)G(x_2(\theta) + \theta) & (\theta \in [2d_0, \infty)). \end{align}$$

**(b)** *Furthermore, $x_0$ and $x_1$ are increasing, while $x_2$ is decreasing, with*

$$\lim_{\theta \to 0^+} x_0(\theta) = -\infty, x_0(d_1) = -d_1, x_0(2d_0) = -d_0;$$

$$\lim_{\theta \to 0^+} x_1(\theta) = d_0, \lim_{\theta \to \infty} x_1(\theta) = d_1;$$

$$x_2(2d_0) = -d_0, \text{ and } \lim_{\theta \to \infty} x_2(\theta) = -d_1.$$

*Proof.* **(a)** First, observe that the nondecreasing properties of $l(\cdot)$ and $u(\cdot)$ established in Lemma 4 imply that $l^{-1}(\cdot)$ and $u^{-1}(\cdot)$ are equally nondecreasing. Along with Lemma 3, it also follows that $l(x)$ varies continuously from 0 to $+\infty$ as $x \in \Re$; and that $u(x)$ varies continuously from $a$ to $+\infty$ as $x \in \Re$. Set

$$\begin{align} x_0(\theta) &= u^{-1}(\theta) - \theta; \quad \theta \in (a, 2d_0]; \\ x_1(\theta) &= l^{-1}(\theta) - \theta; \quad \theta \in [0, \infty); \\ x_2(\theta) &= u^{-1}(\theta) - \theta; \quad \theta \in [2d_0, \infty). \end{align}$$

We hence obtain

$$\begin{align} C(\theta) &= P_\theta(l(X) \leq \theta \leq u(X)) \\ &= P_\theta(u^{-1}(\theta) \leq X \leq l^{-1}(\theta)) \\ &= G(l^{-1}(\theta) - \theta) - G(u^{-1}(\theta) - \theta), \end{align}$$

which is (5) with $u^{-1}(\theta) = -\infty$ for $\theta \in [0, a]$. Furthermore, it must be that:

- $l^{-1}(\theta) - G^{-1}(\frac{1}{2} + \frac{1-\alpha}{2}G(l^{-1}(\theta)) = \theta$, for $\theta \geq 0, \Rightarrow G(l^{-1}(\theta) - \theta) = \frac{1}{2} + \frac{1-\alpha}{2}\,G(l^{-1}(\theta)) \Rightarrow (ii)$;
- $u^{-1}(\theta) - G^{-1}(\alpha G(u^{-1}(\theta)) = \theta$, for $\theta \in (a, 2d_0), \Rightarrow G(u^{-1}(\theta) - \theta) = \alpha G(u^{-1}(\theta)) \Rightarrow (i)$
- $u^{-1}(\theta) + G^{-1}(\frac{1}{2} + \frac{1-\alpha}{2}G(u^{-1}(\theta))) = \theta$, for $\theta \geq 2d_0, \Rightarrow G(u^{-1}(\theta) - \theta) = \frac{1}{2} + \frac{1-\alpha}{2}G(u^{-1}(\theta)) \Rightarrow (iii)$.

**(b)** The right hand sides of (i), (ii), and (iii) increase in $\theta$ given that both $l^{-1}(\theta)$ and $u^{-1}(\theta)$ increase in $\theta$. This implies that, in terms of $\theta$ ($\geq 0$), $x_0(\theta)$, $x_1(\theta)$ increase, while $x_2(\theta)$ decreases. Now, observe that $u(d_0) = 2d_0 \Rightarrow u^{-1}(2d_0) = d_0$, and $x_0(2d_0) (= x_2(2d_0)) = d_0 - 2d_0 = -d_0$. Similarly, $u(0) = d_1 \Rightarrow x_0(d_1) = -d_1$, and $\lim_{x \to -\infty} u(x) = a \Rightarrow \lim_{\theta \to a^+} x_0(\theta) = -\infty$. Now, making use of (ii) and (iii) and the limiting properties $\lim_{\theta \to \infty} l^{-1}(\theta) = \infty$, $\lim_{\theta \to \infty} u^{-1}(\theta) = \infty$, we infer that $\lim_{\theta \to \infty} x_1(\theta) = d_1$, and $\lim_{\theta \to \infty} x_2(\theta) = -d_1$. Finally, the property $\lim_{\theta \to 0^+} x_1(\theta) = d_0$ follows directly from (ii), or again from part (a) of Lemma 3. □

*Remark* 1. Above, only the properties relative to $x_0$, as well as those of the coverage $C(\theta)$ for $\theta \leq 2d_0$ require the logconcavity of $g$. Analogously, some results below (e.g., Corollary 1) do not require the additional assumption of logconcavity.

*Remark* 2. From (5) and the properties of $x_0, x_1$, and $x_2$ of Lemma 5, it follows that $C(\cdot)$ is a continuous function on $[0, \infty)$, with the exception of a discontinuity at $\theta = a$, and when $a > 0$. In this case, we have $\lim_{\theta \to a^-} x_0(\theta) = -\infty$, and $\lim_{\theta \to a^+} x_0(\theta) = u^{-1}(a) - a$, which will lead to a drop of $G(u^{-1}(a) - a)$ in the coverage at $\theta = a$. An interesting example of such a discontinuity occurs for a Laplace model (see Example 3, Figure 2, and Remark 3).

*Remark* 3. It is pertinent here to discuss the behaviour of $u(x)$ as $x \to -\infty$. In particular, we wish to single out cases where $a > 0$, which will imply that $u(x) \geq \theta$ for any $\theta \in [0, a]$ (i.e., $I_{\pi^*}(x)$ does not underestimate $\theta$ for such $\theta$'s, and coverage will occur as soon as underestimation does not occur; see (5)). As an example, consider a Laplace model with $g(z) = \bar{G}(z) = \frac{1}{2} e^{-z}$; for $z > 0$; and which leads to $u(x) = -\ln(\alpha)$ for all $x < 0$, hence $a = -\ln(\alpha) > 0$. Part (a) of Lemma 3 provides a way to verify this directly. Alternatively, observe that the posterior survivor function of $\theta$ is given by

$$P(\theta \geq y | x) = \frac{\bar{G}(y-x)}{\bar{G}(-x)} = e^{-y}; \text{ for } y > 0, x < 0.$$

Thus, for $x < 0$, the posterior distribution does not vary and yields a constant $I_{\pi^*}(x) = [0, -\ln(\alpha)]$.

Analogously, logconcave densities $g$ with exponential tails will lead to a similar non-zero limit at $-\infty$. A family of such densities, which will lead to $u(x) \to -\ln(\alpha)$ as $x \to -\infty$, is given by $g(z) = P(|z|)e^{-|z|}$; with $P(\cdot)$ nondecreasing and logconcave on $(0, \infty)$, $P'(0^+) < P(0)$, and $\frac{P'(z)}{P(z)} \to 0$ as $z \to \infty$. This may be verified by showing that the conditions on $P$ force the density $g$ to be logconcave, and that the posterior survivor function $P(\theta \geq y|x)$ converges (for $y > 0$) as above to $e^{-y}$ when $x \to -\infty$. A simple example is given by $P(z) = \frac{1}{4}(|z| + 1)$, that is $g(z) = \frac{1}{4}(|z| + 1)e^{-|z|}$.

On the other hand, if $h'$ is unbounded where $h \equiv -\ln(g)$, then $u(x) \to 0$ as $x \to -\infty$. To prove this, it is sufficient to show that $P(\theta \geq y|x) \to 0$ as $x \to -\infty$ for all $y > 0$. But notice that (for $y > 0$)

$$\lim_{x \to -\infty} P(\theta \geq y|x) = \lim_{x \to -\infty} \frac{\bar{G}(y-x)}{\bar{G}(-x)} = \lim_{x \to \infty} \frac{g(y+x)}{g(x)} = 0,$$

given that $\lim_{x \to \infty} \{h(y+x) - h(x)\} \geq \lim_{x \to \infty} y h'(x) = \infty$.

Finally, we emphasize that the assumption of logconcavity is indeed required for the upper bound $u(\cdot)$ to increase on $(-\infty, d_0)$. As for interesting counterexamples, we point out that $u(x) \to \infty$ as $x \to -\infty$ for non-logconcave densities $g$ such that $\lim_{x \to \infty} \frac{g(y+x)}{g(x)} = 1$, for all $y > 0$. This is the case for instance of Student densities with $\nu \geq 1$ degrees of freedom, as remarked upon previously by [10].

By virtue of Lemma 5, we now obtain a first set of new results for the coverage $C(\theta)$, $\theta \geq 0$.

**Corollary 1.** *For $X \sim g(x - \theta)$, $\theta \geq 0$, $g$ unimodal, symmetric about $0$, the frequentist coverage $C(\theta)$ of the HPD credible set $I_{\pi^*}(X)$ satisfies the properties:*

(a) $C(\theta)$ *increases in* $\theta$ *for* $\theta \geq 2d_0$;
(b) $C(\theta) \leq 1 - \frac{\alpha}{2}$ *for all* $\theta \geq 0$;
(c) $C(\theta) \geq \frac{1}{1+\alpha}$ *for* $\theta \in [0, a]$.

*Proof.* **(a)** Immediate given Lemma 5's representation $G(x_1(\theta)) - G(x_2(\theta))$ for the coverage, and given that $G(\cdot)$ and $x_1(\cdot)$ are increasing on $\mathbb{R}$, while $x_2(\cdot)$ is decreasing on $[2d_0, \infty)$.

**(b)** It follows from Lemma 5 that $C(\theta) \leq G(x_1(\theta)) \leq G(d_1) = 1 - \frac{\alpha}{2}$, for all $\theta \geq 0$.

**(c)** It follows from Lemma 5 that $C(\theta) = G(x_1(\theta)) \geq G(x_1(0)) = G(d_0) = \frac{1}{1+\alpha}$, for all $\theta \in [0, a]$. □

**Corollary 2.** *For $X \sim g(x - \theta)$, $\theta \geq 0$, $g$ unimodal, symmetric about $0$, and logconcave, the frequentist coverage $C(\theta)$ of the HPD credible set $I_{\pi^*}(X)$ satisfies the properties:*

(a) $C(\theta) \leq 1 - \alpha$ *for all* $\theta \geq d_1$;
(b) $C(\theta) \geq 1 - \alpha$ *for all* $\theta \leq G^{-1}(1 - \frac{\alpha^2}{1+\alpha}) - d_0 = d_2$ *(say)*.

*Proof.* **(a)** It suffices to show that $C(\theta) \leq 1 - \alpha$ for $\theta \in [d_1, 2d_0)$ since $C(\theta) \leq 1 - \alpha$ for all $\theta \geq 2d_0$, given part (d) of Lemma 3 and part (a) of Corollary 1. Observe that Lemma 5's increasing property of $x_0$ implies $d_1 \leq 2d_0$ given that $x_0(d_1) = -d_1 \leq -d_0 = x_0(2d_0)$. Furthermore, the properties of Lemma 5 tell us that: $x_1(\theta) \leq d_1$ and $x_0(\theta) \geq -d_1$ for $\theta \in [d_1, 2d_0]$. Finally as a consequence, we obtain from (5) with $\theta \in [d_1, 2d_0) \Rightarrow \theta \geq d_1 = u(0) \geq a$:

$$C(\theta) = G(x_1(\theta)) - G(x_0(\theta)) \leq G(d_1) - G(-d_1) = 1 - \alpha, \text{ for all } \theta \in [d_1, 2d_0).$$

**(b)** Take $\theta \geq a$ as the result is already established for $\theta \leq a$, given that $1 - \alpha \leq \frac{1}{1+\alpha}$. Now define $\theta_2$ such that $G(x_0(\theta_2)) = \frac{\alpha^2}{1+\alpha}$. It is the case that $G(x_0(\theta)) \leq \frac{\alpha^2}{1+\alpha}$ for $\theta \leq \theta_2$, given that $x_0(\cdot)$ is a nondecreasing function. Furthermore, for $\theta \leq \theta_2$, we have

$$\begin{aligned} C(\theta) &\geq G(x_1(\theta)) - G(x_0(\theta)) \\ &\geq G(x_1(0)) - \frac{\alpha^2}{1+\alpha} = G(d_0) - \frac{\alpha^2}{1+\alpha} = \frac{1}{1+\alpha} - \frac{\alpha^2}{1+\alpha} = 1 - \alpha. \end{aligned}$$

There remains to show that $\theta_2 = d_2$. But this is the case by definitions of $x_0$ and $d_2$, with $G^{-1}(\Delta) = -G^{-1}(1-\Delta)$ for all $\Delta \in (0,1)$, since

$$\begin{aligned}\frac{\alpha^2}{1+\alpha} &= G(x_0(\theta_2)) = \alpha G(x_0(\theta_2) + \theta_2) \Rightarrow G^{-1}\left(\frac{\alpha}{1+\alpha}\right) \\ &= x_0(\theta_2) + \theta_2 \Rightarrow \theta_2 = -d_0 - G^{-1}\left(\frac{\alpha^2}{1+\alpha}\right) = d_2.\end{aligned}$$

□

*Remark* 4. Observe that the analysis within the proof of part (b) of the previous corollary tells us also that

$$\sup_{\theta \geq 0} C(\theta) \geq C(d_2) \geq G(x_1(d_2)) - \frac{\alpha^2}{1+\alpha}. \qquad (6)$$

**Example 1.** For sake of illustration of the results obtained up to now, take a normal model with $X|\theta \sim N(\theta, 1)$, and $1-\alpha = 0.90$. Evaluations give us $d_1 = G^{-1}(0.95) \approx 1.645$, $2d_0 = G^{-1}(\frac{1}{1.1}) \approx 2.68$, and $d_2 = G^{-1}(1 - \frac{0.01}{1.1}) - d_0 \approx 1.03$. Results from Corollaries 1 and 2 thus imply that the frequentist coverage $C(\theta)$ necessarily exceeds the nominal coverage 0.90 for $\theta \leq d_2 \approx 1.03$; i.e., for values of the mean $\theta$ that are within a little more than one standard deviation of the lower bound of the parameter space. On the other hand, the exact coverage falls below the nominal coverage and increases for $\theta \geq d_1 \approx 1.645$. Furthermore, we may infer that: **(i)** $C(0) = 0.\overline{90}$, **(ii)** $\sup_{\theta \geq 0} C(\theta) \leq 0.95$, and from (6) **(iii)** $\sup_{\theta \geq 0} C(\theta) \geq G(x_1(d_2)) - \frac{0.01}{1.1} \approx G(x_1(1.03)) - \frac{0.01}{1.1} \approx 0.939$. These above features, along with the increasing property of $C(\theta)$ for $\theta \in [2d_0, \infty)$ (Corollary 1, part a), are illustrated in Figure 1.

The next few results, culminating in Corollary 3, show that the minimum coverage is attained at $\theta = 2d_0$, and is bounded between $1 - \frac{3\alpha}{2}$ and $1 - \frac{3\alpha}{2} + \frac{\alpha^2}{1+\alpha}$. The first of these results gives further information about the behaviour of $C(\theta)$ for $\theta \leq d_1$.

**Lemma 6.** *For $X \sim g(x - \theta)$, $\theta \geq 0$, $g$ unimodal, symmetric about 0, and logconcave, the frequentist coverage $C(\theta)$ of the HPD credible set $I_{\pi^*}(X)$ satisfies the property:*

$$C(\theta) \geq C(2d_0) \text{ for all } \theta \leq d_1.$$

*Proof.* Take $\theta \in [0, d_1]$. From (5), we infer that $C(\theta) \geq C(2d_0)$ as long as

$$G(x_1(\theta)) + G(x_0(2d_0)) \geq G(x_1(2d_0)) + G(x_0(\theta)). \qquad (7)$$

□

Lemma 5's properties of $x_0$ and $x_1$ imply that $G(x_1(\theta)) + G(x_0(2d_0)) \geq G(d_0) + G(-d_0) = 1$, and $G(x_1(2d_0)) + G(x_0(\theta)) \leq G(d_1) + G(-d_1) = 1$; which implies (7) and our desired result.

**Lemma 7.** *For $X \sim g(x - \theta)$, $\theta \geq 0$, $g$ unimodal, symmetric about $0$, and logconcave, the frequentist coverage $C(\theta)$ of the HPD credible set $I_{\pi^*}(X)$ satisfies the properties:*

(i)
$$C(\theta) \geq C(2d_0) \geq 1 - \frac{3\alpha}{2}, \text{ for all } \theta \geq 2d_0.$$

(ii)
$$C(2d_0) \leq 1 - \frac{3\alpha}{2} + \frac{\alpha^2}{1+\alpha}.$$

*Proof of (i).* Since $C(\cdot)$ increases on $[2d_0, \infty)$ (Corollary 1, part a), it suffices to show that
$$C(2d_0) \geq 1 - \frac{3\alpha}{2}. \tag{8}$$

From (5) and the other properties of Lemma 5, we obtain

$$\begin{aligned} C(2d_0) &= G(x_1(2d_0)) - G(x_2(2d_0)) \\ &= \frac{1-\alpha}{2} \{G(x_1(2d_0) + 2d_0) + G(x_2(2d_0) + 2d_0)\} \\ &\geq \frac{1-\alpha}{2} \left\{ G(3d_0 + \frac{1}{1+\alpha}) \right\}. \end{aligned}$$

From this, we obtain that condition (8) is equivalent to $G(3d_0) \geq \frac{1-3\alpha^2}{1-\alpha^2}$, and we conclude the proof by establishing the stronger property (given that $d_0 > 0$)

$$G(2d_0) \geq \frac{1 - 3\alpha^2}{1 - \alpha^2}.$$

With the property $G(d_0) = \frac{1}{1+\alpha}$, rewritten as $\alpha = \frac{\bar{G}(d_0)}{1-\bar{G}(d_0)}$, we infer that

$$G(2d_0) \geq \frac{1-3\alpha^2}{1-\alpha^2} \Leftrightarrow \frac{\bar{G}(2d_0)}{\bar{G}(d_0)} \leq \frac{2\alpha}{1-\alpha} \Leftrightarrow \frac{\bar{G}(2d_0)}{\bar{G}(d_0)} \leq \frac{2\bar{G}(d_0)}{1 - 2\bar{G}(d_0)}.$$

Finally, the proof of (i) is complete since (3) implies the above chain of implications, as well as (8). □

*Proof of (ii).* From (5) and the other properties of Lemma 5, we obtain in a straightforward manner:

$$\begin{aligned} C(2d_0) &= G(x_1(2d_0)) - G(x_0(2d_0)) \\ &\leq G(d_1) - G(-d_0) \\ &= 1 - \frac{\alpha}{2} - \frac{\alpha}{1+\alpha} = 1 - \frac{3\alpha}{2} + \frac{\alpha^2}{1+\alpha}. \end{aligned}$$

The last piece of the analysis consists in showing that the frequentist coverage $C(\theta)$ decreases on $(d_1, 2d_0)$. The proof of the next result relies in part on several lemmas which are stated and proven in the Appendix. □

**Theorem 1.** *For $X \sim g(x - \theta)$, $\theta \geq 0$, $g$ unimodal, symmetric about $0$, and logconcave, the frequentist coverage $C(\theta)$ of the HPD credible set $I_{\pi^*}(X)$ decreases on $(d_1, 2d_0)$.*

*Proof.* Take $\theta \in (d_1, 2d_0)$. Make use (5) and part (c) of Corollary 4 to obtain

$$\begin{aligned} C'(\theta) &\leq x_1'(\theta)\, g(x_1(\theta)) - x_0'(\theta)\, g(x_0(\theta)) \\ &\leq \frac{(1-\alpha)g(d_0)\, g(d_0+d_1)}{2g(d_0) - (1-\alpha)g(d_0+d_1)} - \frac{\alpha}{1-\alpha} g(d_1)\,. \end{aligned}$$

From this, we obtain that the property $C'(\theta) \leq 0$ for $\theta \in (d_1, 2d_0)$ is equivalent to the inequality

$$(1-\alpha)^2\, g(d_0)\, g(d_0+d_1) + \alpha(1-\alpha)g(d_1)\, g(d_0+d_1) \leq 2\alpha\, g(d_0)\, g(d_1).$$

Using the fact that $d_1 \geq d_0 > 0$ and the unimodality of $g$, we infer that the above condition is implied by either the condition

$$\frac{g(d_1)}{g(d_0+d_1)} \geq \frac{1-\alpha}{2\,\alpha};$$

or, given part (b) of Lemma 1, by the stronger condition

$$\frac{g(d_0)}{g(2d_0)} \geq \frac{1-\alpha}{2\,\alpha}\,.$$

Finally, the result follows with this very last inequality being equivalent to (2) given that $\alpha = \frac{\bar{G}(d_0)}{1-\bar{G}(d_0)}$. $\square$

**Corollary 3.** *For $X \sim g(x-\theta)$, $\theta \geq 0$, $g$ unimodal, symmetric about $0$, and logconcave, the frequentist coverage $C(\theta)$ attains its minimum at $\theta = 2d_0$, and is bounded below by $1 - \frac{3\alpha}{2}$.*

*Proof.* The result is a direct consequence of Theorem 1, Lemma 7, and Lemma 6. $\square$

## 3. Examples and final comments

We conclude with some comments and illustrative examples. We also refer to [3] and [6] for further examples and illustrations.

*Remark* 5. The new lower bound $1 - \frac{3\alpha}{2}$ for the frequentist coverage is an improvement over the existing lower bound $\frac{1-\alpha}{1+\alpha}$ for $\alpha < 1/3$, and a significant improvement for relatively smaller $\alpha$. For instance, with a nominal coverage of $1 - \alpha = 0.90$, the lower bounds are $0.85$ and $0.\overline{81}$ respectively. Also, as alluded to in the introduction and as a consequence of Lemma 7, the bound $1 - \frac{3\alpha}{2}$ is, for $\alpha < 1/3$, fairly sharp especially for relatively smaller $\alpha$ given that $\inf_{\theta \geq 0} C(\theta) \leq C(2d_0) \leq 1 - \frac{3\alpha}{2} + \frac{\alpha^2}{1+\alpha}$. For instance with $1 - \alpha = 0.90$, we obtain $0.85 \leq \inf_{\theta \geq 0} C(\theta) \leq 0.85\overline{90}$. There is some numerical evidence to support the applicability of the lower bound $1 - \frac{3\alpha}{2}$ for some models that are not logconcave location models. This is perhaps the case notably for estimating a lower bounded normal mean with unknown variance, (a model associated with a pivotal Student distribution which is not logconcave), based on the reported numerical evaluations in [10].

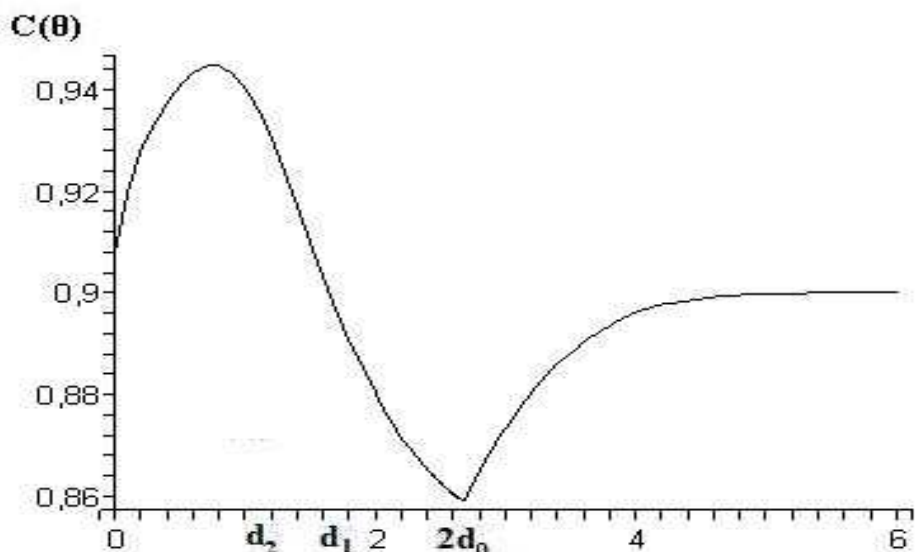

Figure 1 : Frequentist coverage of $I_{\pi^*}(X)$ ; Normal model and $1 - \alpha = 0.90$

**Example 2.** (Normal model; Example 1 continued) As a followup to Example 1, with $X|\theta \sim N(\theta, 1)$, $\theta \geq 0$, and $= 1 - \alpha = 0.90$, the coverage $C(\theta)$ is presented in Figure 1. Notice that the coverage decreases on $(d_1, 2d_0)$ (Theorem 1), and the minimum is attained at $2d_0 \approx 2.68$ and is bounded above by 0.85 (Corollary 3). In fact, the exact minimal coverage $C(2d_0)$ is about equal to 0.859, almost the same as the upper bound given above in Remark 5, or Lemma 7.

**Example 3.** (Laplace model) Consider $X|\theta \sim \text{Laplace}(\theta)$ with $g(z) = \bar{G}(z) = \frac{1}{2}e^{-z}$ for $z > 0$, and $G^{-1}(t) = -\log(2(1-t))$ for $t \geq 1/2$. In view of these closed forms, explicit forms for $I_{\pi^*}(X)$ and its coverage $C(\theta)$ are available. With further details provided in [3] or [6], we obtain for instance from (5), with $a = -\ln(\alpha)$ (see Remark 2), $d_0 = G^{-1}(\frac{1}{1+\alpha}) = \ln(\frac{1+\alpha}{2\alpha})$, $d_1 = G^{-1}(1 - \frac{\alpha}{2}) = -\ln(\alpha) (= a)$:

$$C(\theta) = \begin{cases} 1 - \frac{\alpha e^\theta}{2e^\theta - (1-\alpha)} & \text{if } 0 \leq \theta \leq -\ln(\alpha) \\ 1 - \frac{\alpha e^\theta}{2e^\theta - (1-\alpha)} - \frac{\alpha e^{-\theta}}{2(\alpha - \sqrt{\alpha^2 - \alpha e^{-\theta}})} & \text{if } -\ln(\alpha) < \theta \leq 2\ln(\frac{1+\alpha}{2\alpha}) \\ 1 - \frac{\alpha e^\theta}{2e^\theta - (1-\alpha)} - \frac{(1-\alpha)e^{-\theta}}{2(\sqrt{\alpha^2 + 2(1-\alpha)e^{-\theta}} - \alpha)} & \text{if } \theta \geq 2\ln(\frac{1+\alpha}{2\alpha}) \,. \end{cases}$$

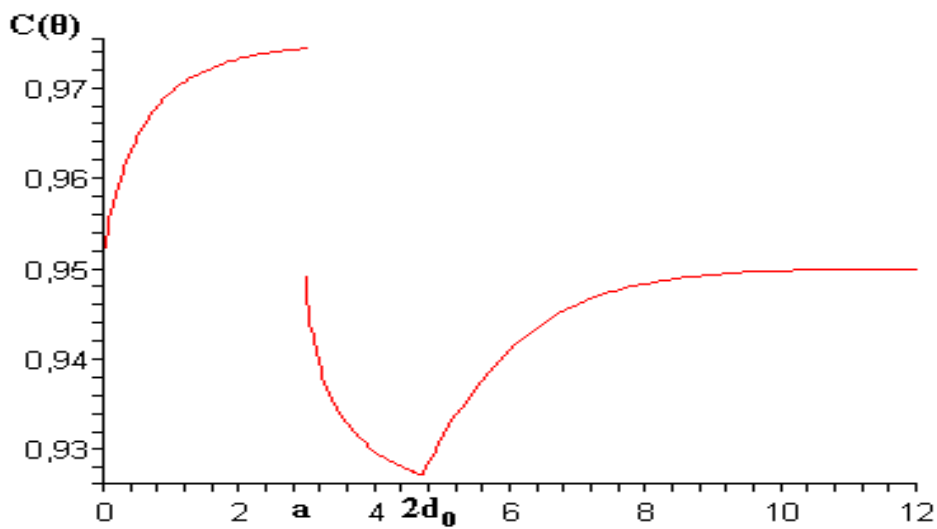

Figure 2 : Frequentist coverage of $I_{\pi^*}(X)$ ;
Laplace model and $1-\alpha = 0.95$

Furthermore, we can show here that $C(\theta)$ is precisely increasing for $\theta \in (0, d_1)$ and $\alpha \leq 1/3$ (we already know it is decreasing on $(d_1, 2d_0)$, and increasing on $(2d_0, \infty)$). This implies here as well that $\sup_{\theta \geq 0} C(\theta) = C(d_1)$. A graph of this coverage is presented in Figure 2 for $1 - \alpha = 0.95$. We recognize several of the established features of $C(\theta)$, such as those illustrated in Examples 1 and 2. Observe the drop at $a = -\ln(0.05) \approx 2.996$ of $G(x_0(a)) = G(x_0(d_1)) = G(-d_1) = \frac{\alpha}{2} = 0.025$ (we made use here of (5) and Lemma 5's properties of $x_0$). Notice also that the minimum is attained at $2d_0 = 2\ln(\frac{1.05}{2(0.05)}) \approx 4.70$ and equal to $C(2d_0) \approx 0.92727$ (using directly the expression above), which is very close to the upper bound $0.92738$ given by Lemma 7. Finally, in this case, $\sup_{\theta \geq 0} C(\theta) = C(-\ln(0.05)) \approx 0.9744$; which of course is less (but only a little less) than Corollary 1's upper bound $0.975$.

Finally, the findings in this paper does provide a sharper description of the frequentist coverage properties of the HPD credible interval $I_\pi^*(X)$ with the improved lower bound on the minimal coverage mitigating in favour of desirable features (of course, added to the fact that the interval $I_\pi^*(X)$ has exact credibility for a uniform prior on $[0, \infty)$). However, one can turn around the argument to point out the non-conformity of the frequentist and nominal coverage of $I_\pi^*(X)$ in the worse case scenario $\theta = 2d_0$. For instance, if $1 - \alpha = 1/3$, Lemma 7 implies that $\inf_{\theta \geq 0} C(\theta)$ is at most $7/12$ (at least $1/2$), in other words a departure of at least $1/12$ between nominal and frequentist coverage at $\theta = 2d_0$.

## 4. Appendix

All the results in this Appendix are called upon in Theorem 1 and are established under the same assumptions as those of Theorem 1, that is: $X \sim g(x-\theta)$, $\theta \geq 0$, $g$ unimodal, symmetric about 0, and logconcave.

**Lemma 8.** *For $\theta \in [d_1, 2d_0]$,*

(a) *the function $\frac{g(x_0(\theta))}{g(x_0(\theta)+\theta)}$ increases in $\theta$;*
(b) *the function $\frac{g(x_1(\theta))}{g(x_1(\theta)+\theta)}$ increases in $\theta$, and is bounded below by $\frac{g(d_0)}{g(d_0+d_1)}$.*

*Proof.* **(a)** Follows from the unimodality of $g$ and the properties of $x_0$ since $x_0(\theta) < 0 \leq x_0(\theta) + \theta$ for $\theta \in [d_1, 2d_0]$.

**(b)** The increasing property follows from Lemma 1 (part c) and Lemma 5's increasing property of $x_1$. With the addition of part (a) of Lemma 1, the lower bound is valid, since for $\theta \in [d_1, 2d_0]$,

$$\frac{g(x_1(\theta))}{g(x_1(\theta) + \theta)} \geq \frac{g(x_1(d_1))}{g(x_1(d_1) + d_1)} \geq \frac{g(d_0)}{g(d_0 + d_1)}.$$

□

**Lemma 9.** *For $\theta \in [d_1, 2d_0]$, we have*

(a) $x_1'(\theta) = \frac{(1-\alpha)g(x_1(\theta)+\theta)}{2g(x_1(\theta)) - (1-\alpha)g(x_1(\theta)+\theta)}$;
(b) $x_0'(\theta) = \frac{\alpha g(x_0(\theta)+\theta)}{g(x_0(\theta)) - \alpha g(x_0(\theta)+\theta)}$.

*Proof.* The expressions follow by differentiation of Lemma 5's implicit equations for $x_0$ and $x_1$. □

**Corollary 4.** *For $\theta \in [d_1, 2d_0]$,*

(a) *both $x_1'(\theta)$ and $x_0'(\theta)$ decrease in $\theta$;*
(b) $x_1'(\theta) \leq \frac{(1-\alpha)g(d_0+d_1)}{2g(d_0)-(1-\alpha)g(d_0+d_1)}$, *and* $x_0'(\theta) \geq \frac{\alpha}{1-\alpha}$.
(c) $x_1'(\theta)\, g(x_1(\theta)) \leq \frac{(1-\alpha)g(d_0)\, g(d_0+d_1)}{2g(d_0)-(1-\alpha)g(d_0+d_1)}$, *and* $x_0'(\theta)\, g(x_0(\theta)) \geq \frac{\alpha}{1-\alpha}\, g(d_1)$.

*Proof.* **(a)** Follows directly from Lemma 8 applied to $\frac{1}{x_1'(\theta)}$ and $\frac{1}{x_0'(\theta)}$ as given in Lemma 9.

**(b)** The upper bound for $x_1'$ follows from Lemma 9 and Lemma 8's lower bound for $\frac{g(x_1(\theta))}{g(x_1(\theta)+\theta)}$. Exploiting the fact that $x_0'$ is decreasing on $[d_1, 2d_0]$, Lemma 9, and the equality $x_0(2d_0) = -d_0$, the lower bound for $x_0'$ is obtained as:

$$x_0'(\theta) \geq x_0'(2d_0) = \frac{\alpha g(x_0(2d_0) + 2d_0)}{g(x_0(2d_0)) - \alpha g(x_0(2d_0) + 2d_0)} = \frac{\alpha g(d_0)}{g(-d_0) - \alpha g(d_0)} = \frac{\alpha}{1-\alpha}.$$

**(c)** Since $x_1(\theta) \in [d_0, d_1]$ (Lemma 5), we have with $g$ being unimodal and $d_0 > 0$: $g(x_1(\theta)) \leq g(d_0)$ for all $\theta \in [d_1, 2d_0]$ (and actually all $\theta \geq 0$ as well). The result then follows directly from (b). □

Similarly, again making use of Lemma 5 and part (b), we have $x_0(\theta) \in [-d_1, -d_0]$ for $\theta \in [d_1, 2d_0]$, and $x'_0(\theta)g(x_0(\theta)) \geq \frac{\alpha}{1-\alpha} g(-d_1) = \frac{\alpha}{1-\alpha} g(d_1)$, for all $\theta \in [d_1, 2d_0]$. □

## Acknowledgements

The comments and suggestions of an anonymous referee are gratefully acknowledged. For the work of Marchand, the support of NSERC of Canada is gratefully acknowledged, while for the work of Strawderman, the support of NSA Grant 03G-1R is gratefully acknowledged.